\vsize=9.0in\voffset=1cm
%\nopagenumbers
\looseness=2

% formato JGR

\message{fonts,}

\font\tenrm=cmr10
\font\ninerm=cmr9
\font\eightrm=cmr8
\font\teni=cmmi10
\font\ninei=cmmi9
\font\eighti=cmmi8
\font\ninesy=cmsy9
\font\tensy=cmsy10
\font\eightsy=cmsy8
\font\tenbf=cmbx10
\font\ninebf=cmbx9
\font\tentt=cmtt10
\font\ninett=cmtt9

\font\ninesl=cmsl9
\font\eightsl=cmsl8

\font\nineit=cmti9
\font\eightit=cmti8
 % for titles

\skewchar\ninei='177 \skewchar\eighti='177
\skewchar\ninesy='60 \skewchar\eightsy='60

\def\eightpoint{\def\rm{\fam0\eightrm} % pasa a font de 8 puntos
\normalbaselineskip=9pt
\normallineskiplimit=-1pt
\normallineskip=0pt

\textfont0=\eightrm \scriptfont0=\sevenrm \scriptscriptfont0=\fiverm
\textfont1=\ninei \scriptfont1=\seveni \scriptscriptfont1=\fivei
\textfont2=\ninesy \scriptfont2=\sevensy \scriptscriptfont2=\fivesy
\textfont3=\tenex \scriptfont3=\tenex \scriptscriptfont3=\tenex
\textfont\itfam=\eightit  \def\it{\fam\itfam\eightit} % \it is family 4
\textfont\slfam=\eightsl \def\sl{\fam\slfam\eightsl} % \sl is family 5

\setbox\strutbox=\hbox{\vrule height6pt depth2pt width0pt}%
\normalbaselines \rm}

\def\ninepoint{\def\rm{\fam0\ninerm} % pasa a font de 9 puntos
\textfont0=\ninerm \scriptfont0=\sevenrm \scriptscriptfont0=\fiverm
\textfont1=\ninei \scriptfont1=\seveni \scriptscriptfont1=\fivei
\textfont2=\ninesy \scriptfont2=\sevensy \scriptscriptfont2=\fivesy
\textfont3=\tenex \scriptfont3=\tenex \scriptscriptfont3=\tenex
\textfont\itfam=\nineit  \def\it{\fam\itfam\nineit} % \it is family 4
\textfont\slfam=\ninesl \def\sl{\fam\slfam\ninesl} % \sl is family 5
\textfont\bffam=\ninebf \scriptfont\bffam=\sevenbf
\scriptscriptfont\bffam=\fivebf \def\bf{\fam\bffam\ninebf} % \bf is family 6
\textfont\ttfam=\ninett \def\tt{\fam\ttfam\ninett} % \tt is family 7

\normalbaselineskip=11pt
\setbox\strutbox=\hbox{\vrule height8pt depth3pt width0pt}%
\let \smc=\sevenrm \let\big=\ninebig \normalbaselines
\parindent=1em
\rm}

\def\tenpoint{\def\rm{\fam0\tenrm} % font de 10 points
\textfont0=\tenrm \scriptfont0=\ninerm \scriptscriptfont0=\fiverm
\textfont1=\teni \scriptfont1=\seveni \scriptscriptfont1=\fivei
\textfont2=\tensy \scriptfont2=\sevensy \scriptscriptfont2=\fivesy
\textfont3=\tenex \scriptfont3=\tenex \scriptscriptfont3=\tenex
\textfont\itfam=\nineit  \def\it{\fam\itfam\nineit} % \it is family 4
\textfont\slfam=\ninesl \def\sl{\fam\slfam\ninesl} % \sl is family 5
\textfont\bffam=\ninebf \scriptfont\bffam=\sevenbf
\scriptscriptfont\bffam=\fivebf \def\bf{\fam\bffam\tenbf} % \bf is family 6
\textfont\ttfam=\tentt \def\tt{\fam\ttfam\tentt} % \tt is family 7

\normalbaselineskip=11pt
\setbox\strutbox=\hbox{\vrule height8pt depth3pt width0pt}%
\let \smc=\sevenrm \let\big=\ninebig \normalbaselines
\parindent=1em
\rm}

\message{fin format jgr}
%end format jgr

\magnification=1200
\font\Bbb=msbm10
\def\C{\hbox{\Bbb C}}
\def\R{\hbox{\Bbb R}}
\def\v{\varphi}
\def\pa{\partial}
\def\b{\backslash}

\setbox3=\vtop{
\hbox{Universit\'e Pierre et Marie Curie, Math\'ematiques, case 247}
\hbox{4 place Jussieu, 75252 Paris, France}
\hbox{e-mail: henkin@math.jussieu.fr,\ michel@math.jussieu.fr}
}

\vskip 2 mm
\centerline{\bf Inverse conductivity problem on
Riemann surfaces}
\vskip 2 mm
\centerline{\bf By Gennadi Henkin and Vincent Michel}

\vskip 2 mm
{\bf Abstract}

An electrical potential $U$ on a bordered Riemann surface $X$ with
conductivity
function $\sigma>0$ satisfies equation $d(\sigma d^cU)=0$. The problem of
effective reconstruction of $\sigma$ from electrical currents measurements
(Dirichlet-to-Neumann mapping) on the boundary:

\noindent
$U\big|_{bX}\mapsto\sigma d^cU\big|_{bX}$ is studied. We extend to the case
of Riemann surfaces the reconstruction scheme given, firstly, by R.Novikov
[N1] for simply connected $X$. We apply for this  new kernels for $\bar\pa$
on the affine algebraic Riemann surfaces constructed in [H2].

\vskip 2 mm
{\bf 0. Introduction}

{\it 0.1. Inverse conductivity problem}.

Let $X$ be bordered oriented real two-dimensional manifold in $\R^3$ equiped
with
a smooth symmetric and positive tensor $\hat\sigma:\ T^*X\to T^*X$ on
cotangent bundle $T^*X$, called anisotropic conductivity tensor on $X$,
$\hat\sigma$ is called  symmetric and positive if
$\hat\sigma a\wedge b=\hat\sigma b\wedge a$ and $\hat\sigma a\wedge a>0$
for any $a,b\in T^*X$.

Let $u$ (correspondingly $U$) be a smooth function on $bX$
(correspondingly on $X$) such that $U\big|_{bX}=u$, called electric
potential on $bX$ and correspondingly on $X$. 1-form $\hat\sigma dU$
 on $X$ is called electrical current on $X$. By Maxwell equation
$$d(\hat\sigma dU)=0\ \ {\rm on}\ \ X.$$
Inverse conductivity problem consists in this case in the following:
what kind of information about $X$ and $\hat\sigma$ can be efficiently
extracted from the knowledge of Dirichlet-to-Neumann mapping
$$u\big|_{bX}\mapsto\hat\sigma dU\big|_{bX}\ \forall\ u\in C^{(1)}(bX).
$$
This important problem in the mathematical setting goes back to the inverse
boundary values problems posed by I.M.Gelfand [G] and by A.P.Calderon [C].

This real problem is deeply related with complex analysis on Riemann
surfaces. The first indication to this relation gives the following
statement, obtaind at first by J.Sylvester [S] for simply connected $X$.

Under conditions above\ $\forall$ couple $(X,\hat\sigma)$ there exist
a unique complex structure $c$ on $X$ and smooth scalar valued, positive
conductivity function $\sigma$ such that the equation
$d(\hat\sigma dU)=0$ takes the form $d(\sigma d^cU)=0$, where
$d^c=i(\bar\pa-\pa)$, $\sigma^2(x)=det\,\hat\sigma(x)$, $x\in X$.

This statement permits to reduce the inverse conductivity problem to the
questions about reconstruction from Dirichlet-to-Neumann mapping of the
genus of $X$, of the complex structure of $X$ and of the scalar conductivity
function $\sigma$ on $X$.

These questions are well answered for the important case when $X$ is a
domain in $\R^2$, due to the sequence of works:
[F1], [F2], [BC1], [SU], [N3], [N1], [GN], [BC2], [N2], [Na].

The exact reconstruction scheme for this case was discovered
by R.Novikov [N1].

Formulated questions are well answered also for the case when
conductivity function $\sigma$ is known to be constant on $X$, i.e.
when  only Riemann surface $X$
must be reconstructed from Dirichlet-to-Neumann data [LU], [Be], [HM].

In this paper we study another important case of this problem, when
bordered two-dimensional manifold $X$ and complex structure on $X$ are
known, but conductivity function $\sigma$ on $X$ must be reconstructed from
Dirichlet-to-Neumann data.

{\it 0.2. Main results}.

We extend here the R.Novikov's reconstruction scheme for the case of
bordered Riemann surfaces. Our method (announced firstly in [H1])
is based on the  appropriate new
kernels for $\bar\pa$  on the affine algebraic Riemann surfaces constructed
in [H2].

By this reason we use some special embedding of $X$ into $\C^2$.

Let $\hat X$ be compactification of $X$ such that
$\hat X=\overline{X\cup X_0}$ be compact Riemann surface of genus $g$.
Let $A=\{A_1,\ldots,A_d\}$ be divisor, generic, effective with support
in $X_0$, consisting of $d=g+2$ points.

By Riemann-Roch formula there exist three independent functions
$f_0,f_1,f_2\in {\cal M}(\hat X)\cap {\cal O}(\hat X\b A)$ having at most
simple poles in the points of $A$.
Without restriction of generality one can put $f_0=const$. Let $V$
be algebraic curve in $\C^2$ of the form
$$V=\{(z_1,z_2)\in\C^2:\ z_1=f_1(x),\ z_2=f_2(x),\ x\in \hat X\b A\}.$$
Let $\tilde V$ be compactification of $V$ in ${\C}P^2$ of the form
$\tilde V=\{w\in {\C}P^2:\ \tilde P(w)=0\}$, where $\tilde P$ is
homogeneous holomorphic polynomial of homogeneous coordinates
$w=(w_0:w_1:w_2)$. Without loss of generality  one can suppose:
functions $f_1$, $f_2$ are such that
\item{  i)} $\tilde V$ intersects ${\C}P^1_{\infty}=\{w\in {\C}P^2:\ w_0=0\}$
 transversally $\tilde V\cap {\C}P^1_{\infty}=\{a_1,\ldots,a_d\}$, where
points $a_j=\bigl(0,1,\lim\limits_{x\to A_j}{f_2(x)\over f_1(x)}\bigr)$,
$j=1,2,\ldots,d$ are different points of ${\C}P^1_{\infty}$.
\item{ ii)} $V=\tilde V\b {\C}P^1_{\infty}$ is connected curve in $\C^2$
with equation $V=\{z\in\C^2:\ P(z)=0\}$, where $P(z)=\tilde P(1,z_1,z_2)$
such that $|{\pa P\over \pa z_1}|\le const(V)|{\pa P\over \pa z_2}|$, if
$|z_1|\ge r_0=const(V)$.
\item{iii)} For any $z^*\in V$, where ${\pa P\over \pa z_2}(z^*)=0$
we have ${\pa^2 P\over \pa z_2^2}(z^*)\ne 0$.

With certain restriction of generality we suppose, in addition, that
\item{ iv)} curve $V$  is a regular curve, i.e.
$grad\,P(z)\ne 0\ \forall\ z\in V$.
This restriction must be eliminated in other publication.

\vskip 2 mm
Let us equip $V$ by euclidean volume form $d d^c|z|^2\big|_V$.

Let $\v\mapsto f=\hat R\v$ be operator for solution of $\bar\pa f=\v$ on $V$

\noindent
from [H2], Proposition 2,
$f\mapsto u=R_{\lambda}f$ be operator for solution of
$(\pa+\lambda dz_1)u=f-{\cal H}f$ on $V$, where ${\cal H}f$ is projection
of $f$ on subspace of holomorphic (1,0)-forms on $\tilde V$
from [H2], Proposition 3, $\v\in L^{\infty}_{1,1}\cap L^1_{1,1}(V)$,
$f\in W_{1,0}^{1,\tilde p}(V)$, $u\in W^{2,\tilde p}(V)$, $\tilde p>2$.

Let
$g_{\lambda}(z,\xi)$, $z,\xi\in V$, $\lambda\in\C$ be kernel of operator
$R_{\lambda}\circ\hat R$ from [H2].

Let $V_X=\{(z_1,z_2)\in V:\ z_1=f_1(x),\ z_2=f_2(x),\ x\in X\}$.

Let $\sigma$ be conductivity function  on $V$ with conditions
$\sigma\in C^{(3)}(V)$, $\sigma>0$ on $V$, $\sigma(z(x))=\sigma(x)$,
$x\in X$,  $\sigma=const$ on $V\b V_X$.

Function $\psi(z,\lambda)$, $z\in V$, $\lambda\in\C$ will be called
Faddeev type function on $V\times\C$ if
$d d^c\psi={d d^c\sqrt{\sigma}\over \sqrt{\sigma}}\psi$ on $V$ and\
$\forall\lambda\in\C$
$e^{-\lambda z}\psi(z,\lambda)\buildrel \rm def \over =\mu(z,\lambda)\to 1$,
$z\to\infty$, $|\bar\pa\mu|=O({1\over {|z|+1}})$, $z\in V$.

\vskip 2 mm
{\bf Theorem}
Under formulated conditions

\item{ I.} There exists unique  Faddeev type function $\psi(z,\lambda)$,
$z\in V$, $\lambda\in\C$.
\item{II.} Function $\psi$ and as a conseqence
conductivity function $\sigma$ can be reconstructed through
Dirichlet-to-Neumann data by the following procedure.

${II}_a$. From Dirichlet-to-Neumann data on $bX$  by
Proposition 3.1 (section 3) one can find restriction
$\psi\big|_{bV_X}$ of the Faddeev type function $\psi(z,\lambda)$,
$z\in V$, $\lambda\in\C$ as a unique solution of the Fredholm integral
equation
$$\psi(z,\lambda)\big|_{bV_X}=e^{\lambda z_1}+
\int\limits_{\xi\in bV_X}e^{\lambda(z_1-\xi_1)}g_{\lambda}(z,\xi)\cdot
(\hat\Phi\psi(\xi)-\hat\Phi_0\psi(\xi)),$$
where $\hat\Phi\psi=\bar\pa\psi\big|_{bV_X}$,
$\hat\Phi_0\psi=\bar\pa\psi_0\big|_{bV_X}$, $d d^c\psi_0\big|_{V_X}=0$,
$\psi_0\big|_{bV_X}=\psi\big|_{bV_X}$.

${II}_b$. Using values of $\psi(z,\lambda)$ in arbitrary point $z^*\in bV_X$
by Proposition 2.2 (section 2) one can find "$\bar\pa$ scattering data":
$$b(\lambda)\buildrel \rm def \over =
\lim\limits_{\scriptstyle z\to\infty
\atop \scriptstyle z\in V}{\bar z_1\over \bar\lambda}
e^{-\bar\lambda \bar z_1}{\pa\psi\over \pa\bar z_1}(z,\lambda)=
(\bar\psi(z^*,\lambda))^{-1}
{\pa\psi\over \pa\bar\lambda}(z^*,\lambda),\ \ z^*\in bV_X$$
with estimate (2.12).

${II}_c$. Using $b(\lambda)$, $\lambda\in\C$ by Proposition 2.3 (section 2)
one can find values of
$$\mu(z,\lambda)\big|_{V_X}=\psi(z,\lambda)e^{-\lambda z_1}\big|_{V_X},\ \
\lambda\in\C$$
as a unique solution of Fredholm integral equation
$$\mu(z,\lambda)=1-{1\over 2\pi i}\int\limits_{\xi\in\C}b(\xi)
e^{\bar\xi\bar z_1-\xi z_1}\overline{\mu(z,\xi)}{d\xi\wedge d\bar\xi\over
{\xi-\lambda}}.$$
From  equality $d d^c\psi={d d^c\sqrt{\sigma}\over \sqrt{\sigma}}\psi$ on $X$
 we find finally
${d d^c\sqrt{\sigma}\over \sqrt{\sigma}}$ on $X$.

\vskip 2 mm
{\it Remarks}

For the case   $V=\C$ the reconstraction scheme I-II for potential $q$
in the Schr\"odinger equation  $-\Delta U+qU=EU$  on $X\subset V$ through
the Dirichlet-to-Neumann data on $bX$ was given for the first time by
R.Novikov [N1].
However, in [N1] this scheme was rigorously justified only for the case
when estimates of the  type (2.12) are available, for example,
if for $E\ne 0$ $\|q\|\le const(E)$.
By the additional result of A.Nachman [Na] the estimates of the type (2.12)
are valid also if $E=0$ and $q={\Delta\sqrt{\sigma}\over \sqrt{\sigma}}$,
$\sigma>0$, $\sigma\in C^{(2)}(X)$.

Part ${II}_a$ in the present paper is completely similar to the related
result of [N1] for $V=\C$.

Part ${II}_b$ of this scheme for $V=\C$ is a consequence of works
R.Beals, R.Coifman [BC1], P.Grinevich, S.Novikov [GN] and R.Novikov [N2].

Part ${II}_c$ of this scheme follows from part  ${II}_b$ and the classical
result of I.Vekua [V].

\vskip 2 mm
\centerline{\bf $\S 1$. Faddeev type function on affine algebraic Riemann
surface.}
\centerline{\bf Uniqueness and existence}

Let $V$ be smooth algebraic curve in $\C^2$ defined in introduction,
equiped by euclidean volume form $d\,d^c|z|^2\big|_V$.

Let $V_0=\{z\in V:\ |z_1|\le r_0\}$, where $r_0$ satisfies condition ii) of
introduction.

\vskip 2 mm
{\it Definition}

Let $q$ be (1,1)-form in $C_{1,1}(\tilde V)$ with support of $q$
in $V_0$. For $\lambda\in\C$ function $z\mapsto\psi(z,\lambda)$, $z\in V$
will be called here the Faddeev type function associated with form
(potential) $q$ on $V$ (and zero level of energy $E$) if
$$-d\,d^c\psi+q\psi=0,\ \ z\in V \eqno(1.1)$$
and function
$\mu=e^{-\lambda z_1}\psi$ satisfies the properties:
$$\mu\in C(\tilde V),\ \ \lim_{\scriptstyle z\to\infty \atop\scriptstyle
z\in V}\mu(z,\lambda)=1\ \ {\rm and}\ \ |\bar\pa\mu|=O\bigl({1\over {1+|z|}}
\bigr),\ z\in V.$$
From  [F1], [F2], [N2] it follows that in the case
$V=C=\{z\in\C^2:\ z_2=0\}$, for almost all $\lambda\in\C$ the Faddeev type
function $\psi=e^{\lambda z_1}\mu$ exists, unique and satisfies the Faddeev
type integral equation
$$\eqalignno{
&\mu(z,\lambda)=1+{i\over 2}\int\limits_{\xi\in V}g(z-\xi,\lambda)
\mu(\xi,\lambda)q(\xi)\ \ {\rm where}\ \ \cr
&g(z,\lambda)={i\over 2(2\pi)^2}\int\limits_{w\in\C}{e^{i(w\bar z+\bar w z)}
dw\wedge d\bar w\over w(\bar w-i\lambda)}\cr}$$
is so called the Faddeev-Green function for operator
$\mu\mapsto\bar\pa(\pa+\lambda dz_1)\mu$ on $V=\C$.

Faddeev type functions $z\mapsto\psi(z,\lambda)$ are especially useful for
solutions of inverse scattering or inverse boundary problems for equation
(1.1) when such  functions exist and unique for any $\lambda\in\C$. It was
remarked by P.Grinevich and S.Novikov [GN] (see also [T])  that for some
continuous $q$ with compact support in $\C$ even with arbitrary
small norm there exists the subset of exceptional $\lambda^*$ for which
the Faddeev type  integral equation  is not uniquely solvable.

From R.Novikov's works [N1], [N2] it follows that  the Faddeev
type functions associated with potential $q$ on $V=\C$ and non-zero level
of energy $E$ exist\ $\forall\lambda\in\C$ if $\|q\|\le const(|E|)$.

From R.Beals, R.Coifman works [BC2] develloped by A.Nachman [Na] and
R.Brown, G.Uhlmann [BU] it follows that for any potential $q$ of the form
$q={d\,d^c\sqrt{\sigma}\over \sqrt{\sigma}}$, where $\sigma\in C^2(\C)$,
$\sigma(z)\equiv const$ if $|z|\ge const$,  Faddeev type function
$z\mapsto\psi(z,\lambda)$ exists and unique for any $\lambda\in\C$.

Proposition 1.1 below gives the uniqueness of the Faddeev type function
on affine algebraic Riemann surface $V$ for potential
$q={d\,d^c\sqrt{\sigma}\over \sqrt{\sigma}}$, where $\sigma\in C^2(V)$,
$\sigma=const$ on $V\b V_X \subset V\b V_0$.
The proof will be based on the approach going back to [BC2] in the case
$V=\C$.

\vskip 2 mm
{\bf Proposition 1.1} (Uniqueness)

Let  $\sigma$ be positive function belonging to $C^{(2)}(V)$ such
that $\sigma\equiv const>0$ on
$$V\b V_X\subset V\b V_0=\cup_{j=1}^dV_j,$$
where $\{V_j\}$ are connected components of $V\b V_0$.

Let $\mu\in L^{\infty}(V)$ such that
${\pa\mu\over \pa z_1}\in L^{\tilde p}(V)$ for some $\tilde p>2$ and $\mu$
satisfies equation
$$\bar\pa(\pa+\lambda dz_1)\mu={i\over 2}q \mu\ \ {\rm where}\ \
q={d\,d^c\sqrt{\sigma}\over \sqrt{\sigma}}\eqno(1.2)$$
and for some $j\in\{1,2,\ldots,d\}$ $\mu(z)\to 0$, $z\to\infty$, $z\in V_j$.

Then $\mu\equiv 0$.

\vskip 2 mm
{\it Remark}

Proposition 1 is still valid if to replace the condition
${\pa\mu\over \pa z_1}\in L^{\tilde p}(V)$, $\tilde p>2$ by the weaker
condition   $\pa\mu\in L^{\tilde p}(V)$.

\vskip 2 mm
{\bf Lemma 1.1}

Put $f=e^{\lambda z_1}\mu$, $f_1=\sqrt{\sigma}{\pa f\over \pa z_1}$,
$f_2=\sqrt{\sigma}{\pa f\over \pa\bar z_1}$,
where $\mu$ satisfies conditions of

\noindent
Proposition 1.1. Then
$$\eqalignno{
&d(\sigma d^cf)=0\ \ {\rm on}\ \ V\ \ {\rm and}\ &(1.3)\cr
&{\pa f_1\over \pa\bar z_1}=q_1 f_2,\ \ {\pa f_2\over \pa z_1}=\bar q_1 f_1,
&(1.4)\cr}$$
where $q_1=-{\pa\log\sqrt{\sigma}\over \pa z_1}$.
Besides, $q_1\in L^p(V_0)\ \forall\ p<2$, $q_1=0$ on $V\b V_0$.

\vskip 2 mm
{\it Proof}

The property $q_1\big|_{V\b V_0}=0$ follows from the property
$\sigma\equiv const$ on $V\b V_0$. Put

\noindent
$V_0^{\pm}=\{z\in V_0:\ \pm\big|{\pa P\over \pa z_2}\big|\ge\pm
\big|{\pa P\over \pa z_1}\big|\}$. Put
$\tilde q_1={\pa\log\sqrt{\sigma}\over \pa z_2}$. Then
$q_1\big|_{V_0^+}\in C^{(1)}(V_0^+)$ and
$\tilde q_1\big|_{V_0^-}\in C^{(1)}(V_0^-)$. The identity
$q_1\big|_{V_0^-}=\bigl({\pa z_1\over \pa z_2}\bigr)^{-1}\tilde q_1$
and property iii) imply that

\noindent
$q_1\in L^p(V_0)\ \forall p<2$.
Equation (1.3) for $f$ is equivalent to the  equation (1.2) for
$\mu=e^{-\lambda z_1}f$. Equation (1.3) for $f$ means that
$$\bar\pa F_1+(\pa\ln\sqrt{\sigma})\wedge F_2=0\ \ {\rm and}\ \
\pa F_2+(\pa\ln\sqrt{\sigma})\wedge F_1=0,\eqno(1.5)$$
where $F_1=\sqrt{\sigma}\cdot \pa f$ and $F_2=\sqrt{\sigma}\cdot \bar\pa f$.
Using $z_1$ as a local coordinate on $V$ we obtain from (1.5) the system
(1.4), where $f_1=dz_1\rfloor F_1$ and $f_2=d\bar z_1\rfloor F_2$.

\vskip 2 mm
{\bf Lemma 1.2}

Put $m_1=e^{-\lambda z_1}f_1$ and $m_2=e^{-\lambda z_1}f_2$, where
$f_1,f_2$ are defined in Lemma 4.1. Then system (1.3) is equivalent to system
$${\pa m_1\over \pa\bar z_1}=q_1m_2,\ \
{\pa m_2\over \pa z_1}+\lambda m_2=\bar q_1m_1.\eqno(1.6)$$
Besides,
$$m_1=\sqrt{\sigma}\bigl(\lambda\mu+{\pa\mu\over \pa z_1}\bigr)\ \ {\rm and}\
m_2=\sqrt{\sigma}{\pa\mu\over \pa\bar z_1}.\eqno(1.7)$$

\vskip 2 mm
{\it Proof}

Putting in (1.4) $f_1=e^{\lambda z_1}m_1$ and $f_2=e^{\lambda z_1}m_2$,
we obtain
$$\eqalign{
&{\pa e^{\lambda z_1}m_1\over \pa\bar z_1}=q_1e^{\lambda z_1}m_2
\Longleftrightarrow {\pa m_1\over \pa\bar z_1}=q_1m_2\ \ {\rm and} \cr
&{\pa e^{\lambda z_1}m_2\over \pa z_1}=\bar q_1e^{\lambda z_1}m_1
\Longleftrightarrow {\pa m_2\over \pa z_1}+\lambda z_2=\bar q_1m_1.\cr}$$
Besides,
$$f_1=\sqrt{\sigma}{\pa f\over \pa z_1}=
\sqrt{\sigma}e^{\lambda z_1}
\bigl(\lambda\mu+{\pa\mu\over \pa z_1}\bigr)=e^{\lambda z_1}m_1,
$$
where $m_1=\sqrt{\sigma}\bigl(\lambda\mu+{\pa\mu\over \pa z_1}\bigr)$ and
$$f_2=\sqrt{\sigma}{\pa f\over \pa\bar z_1}=
\sqrt{\sigma}e^{\lambda z_1}
\bigl({\pa\mu\over \pa\bar z_1}\bigr)=e^{\lambda z_1}m_2,$$
where $m_2=\sqrt{\sigma}{\pa\mu\over \pa\bar z_1}$.

\vskip 2 mm
{\bf Lemma 1.3}

Put $u_{\pm}=m_1\pm e^{-\lambda z_1+\bar\lambda\bar z_1}\bar m_2$.
Then system (1.4) is equivalent to the system
$${\pa u_{\pm}\over \pa\bar z_1}=
\pm q_1e^{-\lambda z_1+\bar\lambda\bar z_1}\bar u_{\pm}.$$

\vskip 2 mm
{\it Proof}

From definition of $u_{\pm}$, using Lemma 1.2, we obtain
$$\eqalign{
&{\pa u_{\pm}\over \pa\bar z_1}={\pa m_1\over \pa\bar z_1}\pm
e^{-\lambda z_1+\bar\lambda\bar z_1}
(\bar\lambda\bar m_2-\bar\lambda\bar m_2+q_1\bar m_1)=\cr
&q_1m_2\pm q_1e^{-\lambda z_1+\bar\lambda\bar z_1}\bar m_1=
\pm q_1e^{-\lambda z_1+\bar\lambda\bar z_1}(\bar m_1\pm
e^{\lambda z_1-\bar\lambda\bar z_1}m_2)=\cr
&\pm e^{-\lambda z_1+\bar\lambda\bar z_1}q_1\bar u_{\pm}.\cr}$$

\vskip 2 mm
{\it Proof of Proposition 1.1}

Let $u_{\pm}=m_1\pm e^{-\lambda z_1+\bar\lambda\bar z_1}\bar m_2$. We will
prove that under conditions of Proposition 1.1\ $\forall\lambda\ne 0$
$u_{\pm}\equiv 0$ together with $m_1$ and $m_2$.

From equality $q_1=-{\pa\log\sqrt{\sigma}\over \pa z_1}$ and Lemma 1.3  we
obtain that
$$\bar\pa u_{\pm}=\pm q_1e^{-\lambda z_1+\bar\lambda\bar z_1}
\bar u_{\pm}d\bar z_1\in L_{0,1}^{\tilde p}(V)\cap L_{0,1}^1(V),\ \
\tilde p>2.\eqno(1.8)$$
From (1.6), (1.7) we obtain that $m_1\in L^{\tilde p}(V)\oplus L^{\infty}(V)$
 and
${\pa m_1\over \pa \bar z_1}\in L^p(V)$, $\forall p\ge 1$.

From (1.6), properties $\mu\in L^{\infty}(V)$,
${\pa\mu\over \pa\bar z_1}\in L^{\tilde p}(V)$, $\mu(z)\to 0$, $z\in V_j$,
$z\to\infty$ and

\noindent
from [H2] Corollary 1.1 we deduce that there exists
$$\lim\limits_{\scriptstyle z\to\infty \atop \scriptstyle z\in V}m_1(z)=
\lim\limits_{\scriptstyle z\to\infty \atop \scriptstyle z\in V}
\sqrt{\sigma}\bigl(\lambda\mu+{\pa\mu\over \pa z_1}\bigr)=
\lim\limits_{\scriptstyle z\to\infty \atop \scriptstyle z\in V}
\sqrt{\sigma}\lambda\mu(z)=0.\eqno(1.9)$$
From (1.8), (1.9)
and generalized Liouville theorem from Rodin [R], Theorem 7.1, it follows
that $u_{\pm}=0$.

It means, in particular, that ${\pa\mu\over \pa\bar
z_1}=0$, which together with (1.9) imply $\mu=0$.  Proposition 1.1 is proved.

\vskip 2 mm
{\bf Proposition 1.2} (Existence)

Let conductivity $\sigma$ on $V$ satisfies conditions of Proposition 1.1.
Then\ \ $\forall\lambda\in\C$ there exists the Faddeev type function
$\psi=\mu e^{\lambda z}$ on $V$ associated with potential
$q={d d^c\sqrt{\sigma}\over \sqrt{\sigma}}$, i.e.
$$\eqalign{
&a)\ \bar\pa(\pa+\lambda dz_1)\mu={i\over 2}q\mu,\ \ {\rm where}\ \
q={d d^c\sqrt{\sigma}\over \sqrt{\sigma}},\cr
&b)\ \mu-1\in W^{1,\tilde p}(V)\cap C(\tilde V)\ \ \forall\tilde p>2,\cr
&\mu-1\in C^{(\infty)}(V\b V_0)\ \ {\rm and}\ \
\mu-1=O\bigl({1\over |z_1|}\bigr),\ \ z\in V\b V_0.\cr}$$
Moreover,\ $\forall\tilde p>2$
$$\eqalign{
&c)\ \|\mu-1\|_{L^{\tilde p}(V)}\le const(V,\sigma,\tilde p)
\{\min\bigl({1\over \sqrt{|\lambda|}},{1\over |\lambda|}\bigr)\},\cr
&\|{\pa\mu\over \pa z_1}\|_{L^{\tilde p}(V)}\le const(V,\sigma,\tilde p)
\{\min(\sqrt{|\lambda|},1)\},\ \ \forall\tilde p>2,\cr
&d)\ \forall\lambda\ \exists\ c\in\C\ \ {\rm such\ that}\ \
\bigl({\pa\mu\over \pa\bar\lambda}-c\bigr)\in W^{1,\tilde p}(V),\cr
&e)\ {\rm under\ additional\ assumption}\ \ \sigma\in C^{(3)}(V)\cr
&\|{\pa^2\mu\over \pa z_1^2}\|_{L^{\tilde p}(V)}\le const(V,\sigma)\cdot
\lambda^{1/(1+1/\tilde p)}.\cr}$$

\vskip 2 mm
{\it Remark}

Proposition 1.2 will be proved by approach going back to
L.Faddeev [F1], [F2] and for the case $V=\C$  develloped in [N1], [N2], [Na].

\vskip 2 mm
{\it Proof}

Put $G_{\lambda}=R_{\lambda}\circ\hat R$, where $\hat R$ is operator, defined
by formula (2.4) from Proposition 2 of [H2], and $R_{\lambda}$ operator,
defined by formula (3.1) from Proposition 3 of [H2]. Let function
$\mu(z,\lambda)$ be such that $\forall\lambda$,
$\mu\in W^{1,\tilde p}(V)\oplus const$
with respect to $z\in V$. Then\ $\forall\lambda\ne 0$ we have properties
$q\mu\in C(V)$, $q\mu=0$ on $V\b V_0$, which imply, in particular, that
mapping $\mu\mapsto q\mu$ is compact operator from
$W^{1,\tilde p}(V)\oplus const$ to $L^{\tilde p}(V)\cap C(V)$.
By Lemmas 2.1, 2.2 and Proposition 2 from [H2] we have
$\hat Rq\mu\in W^{1,\tilde p}(V)$. Proposition 3 ii) implies that
$R_{\lambda}\circ\hat R(q\mu)\in W^{1,\tilde p}(V)$. Hence, mapping
$\mu\mapsto\mu-R_{\lambda}\circ\hat R{i\over 2}q\mu$ is a Fredholm
operator on  $W^{1,\tilde p}(V)\oplus const$. By Proposition 1.1 the operator
$I-R_{\lambda}\circ\hat R({i\over 2}q\cdot)$ is invertible.
 Then there is  unique function $\mu$ such that
$(\mu-1)\in W^{1,\tilde p}(V)$,\ $\forall\tilde p>2$  and
$$\mu=1+R_{\lambda}\circ\hat R({i\over 2}q\mu).\eqno(1.10)$$
Let us check now statement a) of Proposition 1.2. From (1.10) and
Proposition 3i) from [H2] we obtain
$$\eqalign{
&(\pa+\lambda dz_1)\mu=\lambda dz_1+(\pa+\lambda dz_1)
R_{\lambda}\circ\hat R({i\over 2}q\mu)=\cr
&\lambda dz_1+{\cal H}(\hat R({i\over 2}q\mu))+\hat R({i\over 2}q\mu).\cr}
\eqno(1.11)$$
From (1.11) and Proposition 2 of [H2] we obtain
$$\bar\pa(\pa+\lambda dz_1)\mu={i\over 2}q\mu+
\bar\pa(\lambda dz_1+{\cal H}(\hat R({i\over 2}q\mu))={i\over 2}q\mu,$$
where we have used that
${\cal H}(\hat R({i\over 2}q\mu))\in H_{1,0}(\tilde V)$.

Property a) is proved.

For proving of properties b), c) it is sufficient to remark that by
Proposition 3ii) and 3iii) from [H2] we have property c) and
$$\|(1+|z_1|)(\mu-1)\|_{L^{\infty}(V)}\le
const(V,\sigma){1\over \sqrt{|\lambda|}}.$$
To prove property d) let us differentiate (1.10) with respect to
$\bar\lambda$. We obtain
$$\eqalign{
&{\pa\mu\over \pa\bar\lambda}-
R_{\lambda}\circ\bigl(\hat R({i\over 2}q{\pa\mu\over \pa\bar\lambda})\bigr)=
\bar z_1
\bigl(R_{\lambda}\circ\bigl(\hat R({i\over 2}q\mu)\bigr)\bigr)-
R_{\lambda}\bigl(\bar\xi_1\hat R({i\over 2}q\mu)\bigr)=\cr
&\bar z_1(\mu-1)-R_{\lambda}\bigl(\bar\xi_1\hat R({i\over 2}q\mu)\bigr).\cr}
$$
From Proposition 2 of [H2] we deduce that
$$\bar\xi_1\hat R({i\over 2}q\mu)\in
W^{1,\tilde p}_{1,0}(V)\oplus (const)dz_1.$$
From Proposition 3ii) of [H2] and Remark 1 to it we obtain
$$R_{\lambda}\bigl(\bar\xi_1\hat R({i\over 2}q\mu)\bigr)\in
W^{1,\tilde p}_{1,0}(V)\oplus const.$$
From Proposition 1.2 b)   we deduce also
$$\bar z_1(\mu-1)\in W^{1,\tilde p}(V)\oplus const.$$
Hence,
$${\pa\mu\over \pa\bar\lambda}=
\bigl(I-R_{\lambda}\circ\hat R({i\over 2}q\cdot)\bigr)^{-1}
\bigl(\bar z_1(\mu-1)-R_{\lambda}(\bar\xi_1\hat R({i\over 2}q\mu)\bigr)\in
W^{1,\tilde p}(V)\oplus (const).$$
Property e) follows (under condition  $\sigma\in C^{(3)}(V)$) from
Proposition 3iv) of [H2].

\vskip 2 mm
{\bf $\S 2$. Equation ${\pa\mu(z,\lambda)\over \pa\bar\lambda}=b(\lambda)
e^{\bar\lambda\bar z_1-\lambda z_1}\bar\mu(z,\lambda),\ \lambda\in\C$}

For further results it is important to obtain asymptotic development of
${\pa\mu\over \pa\bar z_1}(z,\lambda)$ for $z_1\to\infty$.

\vskip 2 mm
{\bf Proposition 2.1}

Let conductivity function $\sigma$ on $V$ satisfy the conditions of
Proposition 1.1. Let function $\psi(z,\lambda)=\mu(z,\lambda)e^{\lambda z}$
be the Faddeev type function on $V$ associated with potential

\noindent
$q={d d^c\sqrt{\sigma}\over \sqrt{\sigma}}$. Then\ $\forall\lambda\in\C$
function ${\pa\mu(z,\lambda)\over \pa\bar z_1}$ has the following asymptotic
for $z_1\to\infty$:
$$e^{\lambda z_1-\bar\lambda\bar z_1}{\pa\mu\over \pa\bar z_1}\big|_{V_j}=
{B_{1,j}(\lambda)\over \bar z_1}+\sum\limits_{k=2}^{\infty}
{B_{k,j}(\lambda)\over \bar z_1^k},\eqno(2.1)$$
where $|B_{1,j}(\lambda)|=O(\min(1,\sqrt{|\lambda|}))$, $j=1,2,\ldots,d$ and

under additional condition $\sigma\in C^{(3)}(V)$ we have
$$e^{\lambda z_1-\bar\lambda\bar z_1}{\pa^2\mu\over \pa\bar z_1^2}\big|_{V_j}=
\bar\lambda{B_{1,j}\over \bar z_1}+\sum\limits_{k=2}^{\infty}
{{\bar\lambda B_{k,j}(\lambda)-(k-1)B_{k-1,j}(\lambda)}
\over \bar z_1^k},\eqno(2.2)$$
where
$|B_{1,j}(\lambda)|=O(\min(\sqrt{|\lambda|},|\lambda|^{-1/(1+\tilde p)}))$,
$j=1,2,\ldots,d$\ $\forall\tilde p>2$.

For proving of Proposition 2.1 we need the following decomposition
statement for the Faddeev type function $\mu=\psi e^{-\lambda z}$ on
$V\b V_0$.

\vskip 2 mm
{\bf Lemma 2.1}

\item{  i)} Let $\mu$ be function on $V\b V_0$, which satisfies equation
$$\bar\pa(\pa+\lambda dz_1)\mu=0\ \ {\rm on}\ \ V\b V_0 \eqno(2.3)$$
and the property
$$(\mu-1)\big|_{V\b V_0}\in W^{1,\tilde p}(V\b V_0),\
\forall\tilde p>2.$$
Then
$$\eqalign{
&A\buildrel \rm def \over =
{\pa\mu\over \pa z_1}+\lambda\mu\in{\cal O}(\tilde V\b V_0);\ \
A\big|_{V_j}=\lambda+\sum\limits_{k=1}^{\infty}A_{k,j}{1\over z_1^k},\cr
&B\buildrel \rm def \over =e^{-\lambda z_1+\bar\lambda\bar z_1}
{\overline{\pa\mu}\over \pa\bar z_1}={\cal O}(\tilde V\b V_0);\ \
\bar B\big|_{V_j}=\sum\limits_{k=1}^{\infty}B_{k,j}{1\over \bar z_1^k},\cr}
\eqno(2.4)$$
$|z_1|>r_0$, $j=1,2,\ldots,d$.
\item{ ii)} Let
$$M\big|_{V_j}=1+\sum\limits_{k=1}^{\infty}{a_{k,j}\over z_1^k}\ \ {\rm and}\
\ \bar N\big|_{V_j}=1+\sum\limits_{k=1}^{\infty}{b_{k,j}\over \bar z_1^k}$$
be formal series with coefficients $a_{k,j}$ and $b_{k,j}$ determined by
relations
$$\eqalign{
&\lambda a_{k,j}-(k-1)a_{k-1,j}=A_{k,j};\cr
&\bar\lambda b_{k,j}-(k-1)b_{k-1,j}=B_{k,j};\
j=1,2\ldots,d;\ k=1,2,\ldots\cr}$$
Then function $\mu$ has asymptotic decomposition
$$\eqalignno{
&\mu\big|_{V_j}=M\big|_{V_j}+
e^{\bar\lambda\bar z_1-\lambda z_1}\bar N\big|_{V_j},\ \ z_1\to\infty,\ \
{\rm i.e.}&(2.5)\cr
&\mu\big|_{V_j}=M_{\nu}\big|_{V_j}+
e^{\bar\lambda\bar z_1-\lambda z_1}\bar N_{\nu}\big|_{V_j}+
O\bigl({1\over |z_1|^{\nu+1}}\bigr),\ \ {\rm where}\cr
&M_{\nu}\big|_{V_j}=1+\sum\limits_{k=1}^{\nu}{a_{k,j}\over z_1^k},\ \
\bar N_{\nu}\big|_{V_j}=\sum\limits_{k=1}^{\nu}{b_{k,j}\over \bar z_1^k}.
\cr}$$
\item{iii)} Moreover, for any $A\in {\cal O}(\tilde V\b V_0)$,
$A(z)\to\lambda$, $z\to\infty$, there exist
$B\in {\cal O}(\tilde V\b V_0)$,\
$B(z)\to 0$, $z\to\infty$, and $\mu$ such that $\mu$ satisfies (2.3),
(2.4) and $(\mu-1)\in W^{1,\tilde p}(V\b V_0)$.

\vskip 2 mm
{\it Proof}

\item{  i)} From equation (2.3) it follows that
$$\pa\bar\pa(e^{\lambda z_1}\mu(z,\lambda))\big|_{V\b V_0}=0.$$
It means that
$\bar\pa(e^{\lambda z_1}\mu(z,\lambda))=e^{\lambda z_1}\bar\pa\mu$ is
antiholomorphic form  on $V\b V_0$ and

\noindent
$(\pa\mu+\lambda\mu dz_1)$ is
holomorphic form on $V\b V_0$. From this, condition
$\bar\pa\mu\in L_{0,1}^{\tilde p}(V\b V_0)$ and the Cauchy theorem it
follows that
$$\eqalign{
&e^{\lambda z_1}\bar\pa\mu\big|_{V_j}=
e^{\bar\lambda\bar z_1}\bar Bd\bar z_1\big|_{V_j}=
e^{\bar\lambda\bar z_1}
\sum_{k=1}^{\infty}{B_{k,j}\over \bar z_1^k}d\bar z_1\big|_{V_j}\cr
&{\rm and}\ \ (\pa\mu+\lambda\mu dz_1)\big|_{V_j}=A dz_1\big|_{V_j}=
\bigl(\lambda+
\sum_{k=1}^{\infty}{A_{k,j}\over z_1^k}\bigr)d z_1\big|_{V_j}.\cr}$$
\item{ ii)} From (2.4), (2.5) we obtain, at first, that
$$\bar\pa\mu=e^{-\lambda z_1}\bar\pa(e^{\bar\lambda\bar z_1}\bar N_{\nu})+
O\bigl({1\over |\bar z_1|^{\nu+1}}\bigr)$$
and then
$$\mu=M_{\nu}+e^{\bar\lambda\bar z_1-\lambda z_1}\bar N_{\nu}+
O\bigl({1\over |z_1|^{\nu+1}}\bigr).\eqno(2.6)$$
\item{iii)} Let
$$A\big|_{V_j}=\lambda+\sum_{k=1}^{\infty}{A_{k,j}\over z_1^k}\in
L^{\tilde p}(V\b V_0).$$
Proposition 3 iii) from [H2],
applied to (1,0) forms on the complex plan with coordinate $z_1$, gives
that
\ $\forall\lambda\ne 0$ there exists
$\mu\in W^{1,\tilde p}(V\b V_0)\oplus 1$ such that
$A={\pa\mu\over \pa z_1}+\lambda\mu$. Put
$\bar B=e^{\lambda z_1-\bar\lambda\bar z_1}{\pa\mu\over \pa\bar z_1}$. Then
$\bar B\in L^{\tilde p}(V\b V_0)$ and
$${\pa\bar B\over \pa z_1}=
e^{\lambda z_1-\bar\lambda\bar z_1}\bigl({\pa\over \pa\bar z_1}
\bigl(\lambda\mu+{\pa\mu\over \pa z_1}\bigr)\bigr)=
e^{\lambda z_1-\bar\lambda\bar z_1}\bigl({\pa A\over \pa\bar z_1}\bigr)=0,$$
i.e. $B\in {\cal O}(\tilde V\b V_0)$.
By construction $\mu$ satisfies (2.3), (2.4) and $(\mu-1)\in W^{1,\tilde p}$.

\vskip 2 mm
{\bf Lemma 2.2}

Functions $M_{\nu}$ and $N_{\nu}$ from decomposition (2.6) have  the
following properties:
$$\eqalign{
&\forall\ z\in\tilde V\b V_0\ \ \exists\lim\limits_{\nu\to\infty}
\bigl({\pa M_{\nu}\over \pa z_1}+\lambda M_{\nu}\bigr)\buildrel \rm def \over
={\pa M\over \pa z_1}+\lambda M\ \ {\rm and}\cr
&\exists\lim\limits_{\nu\to\infty}
\bigl({\pa N_{\nu}\over \pa z_1}+\lambda N_{\nu}\bigr)\buildrel \rm def \over
={\pa N\over \pa z_1}+\lambda N.\cr}$$
Functions
${\pa M\over \pa z_1}+\lambda M$ and ${\pa N\over \pa z_1}+\lambda N$ belongs to
 to ${\cal O}(\tilde V\b V_0)$ and
$$\eqalignno{
&{\pa\mu\over \pa\bar z_1}=e^{\bar\lambda\bar z_1-\lambda z_1}
({\pa\bar N\over \pa\bar z_1}+\bar\lambda\bar N),&(2.7)\cr
&{\pa\mu \over \pa z_1}+\lambda\mu={\pa M\over \pa z_1}+\lambda M,&(2.8)\cr
&{\pa N \over \pa z_1}+\lambda N\to 0,\ \ {\rm if}\ \ z_1\to\infty.\cr}$$

\vskip 2 mm
{\it Proof}

Let us show that  (2.4) implies (2.7) and (2.8). Indeed,
$$\eqalign{
&{\pa\mu\over \pa\bar z_1}=\lim\limits_{\nu\to\infty}
\bigl(\bar\lambda e^{\bar\lambda\bar z_1-\lambda z_1}\bar N_{\nu}+
e^{\bar\lambda\bar z_1-\lambda z_1}{\pa\bar N_{\nu}\over \pa\bar z_1}\bigr)=
e^{\bar\lambda\bar z_1-\lambda z_1}\lim\limits_{\nu\to\infty}
\bigl({\pa\bar N_{\nu}\over \pa\bar z_1}+\bar\lambda\bar N_{\nu}\bigr)=\cr
&e^{\bar\lambda\bar z_1-\lambda z_1}
\bigl({\pa\bar N\over \pa\bar z_1}+\bar\lambda\bar N\bigr),\cr
&{\pa\mu\over \pa z_1}+\lambda\mu=\lim\limits_{\nu\to\infty}
\biggl({\pa M_{\nu}\over \pa z_1}-\lambda e^{\bar\lambda\bar z_1-\lambda z_1}
\bar N_{\nu}+\lambda M_{\nu}+\lambda
e^{\bar\lambda\bar z_1-\lambda z_1}\bar N_{\nu}\biggr)=\cr
&\lim\limits_{\nu\to\infty}
\bigl({\pa\bar M_{\nu}\over \pa z_1}+\lambda M_{\nu}\bigr)=
{\pa M\over \pa z_1}+\lambda M.\cr}$$
Properties (2.6), (2.7), (2.8), $\bar\pa\mu\in L_{0,1}^{\tilde p}(V)$ and
Riemann extension theorem imply that
${\pa M\over \pa z_1}+\lambda M$ and ${\pa N\over \pa z_1}+\lambda N$ belongs
to ${\cal O}(\tilde V\b V_0)$ and
${\pa N \over \pa z_1}+\lambda N\to 0$, if $z_1\to\infty$.

\vskip 2 mm
{\bf Corollary}

In conditions of Lemmas 2.1, 2.2 we have convergence $M_{\nu}\to M$ and
$N_{\nu}\to N$, $\nu\to\infty$, in general, only in the space of formal
series, in spite of that convergence

\noindent
${\pa M_{\nu}\over \pa z_1}+\lambda M_{\nu}\to
{\pa M\over \pa z_1}+\lambda M$
and
${\pa N_{\nu}\over \pa z_1}+\lambda N_{\nu}\to
{\pa N\over \pa z_1}+\lambda N$
take place in the space ${\cal O}(\tilde V\b V_0)$.

\vskip 2 mm
{\it Simple example}

Put $V_0=\{z\in V:\ |z_1|<1\}$, $\lambda=1$,
$A=\lambda+\sum\limits_{k=1}^{\infty}{A_k\over z_1^k}$ with $A_k=1$,
$k=1,2,\ldots$

By Lemma 2.1 there exists $\mu=M+e^{\bar\lambda\bar z_1-\lambda z_1}\bar N$,
which satisfies (2.3), (2.4), where
$M=1+\sum\limits_{k=1}^{\infty}{a_k\over z_1^k}$ with $a_k$ determined by
relations $a_k-(k-1)a_{k-1}=A_k=1$; $k=1,2,\ldots$.
It gives $a_k=(k-1)!\bigl(1+{1\over 2!}+\ldots +{1\over (k-1)!}\bigr)$.
We have $|a_k|^{1/k}\to\infty$, $k\to\infty$, i.e.  radius of convergence
of serie for $M$ is equal to zero, in spite of that $|A_k|^{1/k}=1$,
$k=1,2,\ldots$.

\vskip 2 mm
{\it Proof of Proposition 2.1}

Estimate for ${\pa\mu\over \pa z_1}$ from Proposition 1.2 c) and the Cauchy
theorem applied to antiholomorphic function
$e^{\lambda z_1-\bar\lambda\bar z_1}{\pa\mu\over \pa\bar z_1}$ implies
development (2.1) and estimate
$$|B_{1,j}(\lambda)|=O(\min(1,\sqrt{|\lambda|})).\eqno(2.9)$$
Estimate for ${\pa^2\mu\over \pa\bar z_1^2}$ from Proposition 1.2 e) and
the Cauchy theorem for antiholomorphic in $z_1$ function
$e^{\lambda z_1-\bar\lambda\bar z_1}{\pa^2\mu\over \pa\bar z_1^2}\big|_{V_j}$
imply development (2.2) and estimate
$$|B_{1,j}(\lambda)|=O(\min(\sqrt{|\lambda|},|\lambda|^{-1/(1+\tilde p)}
)).\eqno(2.10)$$

Proposition 2.1 is proved.

The next proposition gives $\bar\pa$-equation on Faddeev function
$\mu(z,\lambda)$ with respect to parameter $\lambda\in\C$. For the case
$V=\C$ this proposition goes back to Beals, Coifmann [BC1], Grinevich,
S.Novikov [GN] and R.Novikov [N2].

\vskip 2 mm
{\bf Proposition 2.2}

Let conductivity function $\sigma$ on $V$ satisfy the conditions of
Proposition 1.1. Let function $\psi(z,\lambda)=\mu(z,\lambda)e^{\lambda z}$
be the Faddeev type function on $V$ constructed in Proposition 1.2 and
associated with potential $q={d d^c\sqrt{\sigma}\over \sqrt{\sigma}}$. Then\
$\forall z\in V$ the following $\bar\pa$-equation with respect to
$\lambda\in\C$ takes place
$${\pa\mu\over \pa\bar\lambda}=b(\lambda)
e^{\bar\lambda\bar z_1-\lambda z_1}\bar\mu,\eqno(2.11)$$
where
$$\eqalign{
&b(\lambda)\buildrel \rm def \over =
\lim\limits_{\scriptstyle z\to\infty
\atop \scriptstyle z\in V}{\bar z_1\over \bar\lambda}
e^{\lambda z_1-\bar\lambda \bar z_1}{\pa\mu\over \pa\bar z_1},\cr
&|b(\lambda)|\le const(V,\sigma)
\{\min\bigl({1\over \sqrt{|\lambda|}},{1\over |\lambda|}\bigr)\}\ \ {\rm if}\
\  \sigma\in C^{(2)}(V)\ \ {\rm and}\cr
&|b(\lambda)|\le const(V,\sigma,\tilde p)
\{\min\bigl({1\over \sqrt{|\lambda|}},
\bigl({1\over |\lambda|}\bigr)^{1+1/\tilde p}\bigr)\}\ \
{\rm if}\ \ \sigma\in C^{(3)}(V),\ \forall\tilde p>2.\cr}\eqno(2.12)$$

\vskip 2 mm
{\it Remark}

The proof below is a generalization of the R.Novikov proof [N2] of the
corresponding statement for the case $V=\C$.

\vskip 2 mm
{\it Proof}

Equation $d d^c\psi=q\psi$ for the Faddeev type function
$\psi=\mu\cdot e^{\lambda z}$  is equivalent to the equation
$\bar\pa(\pa+\lambda dz_1)\mu={i\over 2}q\mu$. Put
$\psi_{\lambda}=\pa\psi/ \pa\bar\lambda$
and $\mu_{\lambda}=\pa\mu/ \pa\bar\lambda$. By Proposition 1.2 d)\
$\forall\lambda\in\C$ function $\mu_{\lambda}$ belongs to
$W^{1,\tilde p}(V)\oplus const(\lambda)$.  Besides, from equation
$d d^c\psi=q\psi$ it follows the equation
$$d d^c\psi_{\lambda}=q\psi_{\lambda}\ \ {\rm on}\ \ V.$$
From Lemmas 2.1, 2.2, Proposition 2.1 and Proposition 1.2 b),c) we obtain
$$\eqalign{
&{\pa\mu\over \pa\bar z_1}\big|_{V_j}=
e^{\bar\lambda\bar z_1-\lambda z_1}{B_{1,j}(\lambda)\over \bar z_1}+
O\bigl({1\over |z_1|^2}\bigr)\ \  {\rm and}\cr
&\biggl({\pa\mu\over \pa z_1}+\lambda\mu\biggr)\bigg|_{V_j}=\lambda+
{A_{1,j}(\lambda)\over z_1}+
O\bigl({1\over |z_1|^2}\bigr).\cr}\eqno(2.13)$$
From (2.4)-(2.6) and (2.13) we deduce that
$$\eqalignno{
&\mu=1+{a_j(\lambda)\over z_1}+e^{\bar\lambda\bar z_1-\lambda z_1}
{b_j(\lambda)\over \bar z_1}+
O\bigl({1\over |z_1|^2}\bigr),&(2.14)\cr
&{\rm where}\ \ \bar\lambda b_j(\lambda)=B_{1,j},\ \
\lambda a_j(\lambda)=A_{1,j},\ \  j=1,2,\ldots,d.&(2.15)\cr}$$
From (2.14) with help of Proposition 1.2 d) we obtain
$$\eqalign{
&\psi=e^{\lambda z_1}\mu=e^{\lambda z_1}
\bigl(1+{a_j(\lambda)\over z_1}+e^{\bar\lambda\bar z_1-\lambda z_1}
{b_j(\lambda)\over \bar z_1}+O\bigl({1\over |z_1|^2}\bigr)\bigr),\cr
&{\rm and}\ \ \psi_{\lambda}={\pa\psi\over \pa\bar\lambda}=
e^{\bar\lambda\bar z_1}\bigl(b_j(\lambda)+O\bigl({1\over |z_1|}\bigr)\bigr),\
\ z\in V_j.\cr}$$
Put $\mu_{\lambda}=e^{-\lambda z_1}\psi_{\lambda}$. We obtain
$\bar\pa(\pa+\lambda dz_1)\mu_{\lambda}=q\mu_{\lambda}$ and
$$\mu_{\lambda}=e^{\bar\lambda\bar z_1-\lambda z_1}
\bigl(b_j+O\bigl({1\over |z_1|}\bigr)\bigr),\ \ z\in V_j.$$

Proposition 1.1 about uniqueness of the Faddeev type function implies
equality (2.11), where
$$b(\lambda)\buildrel \rm def \over =b_1(\lambda)=\ldots =
b_d(\lambda).\eqno(2.16)$$
We have also equalities
$$B_{1,1}(\lambda)=\ldots =B_{1,d}(\lambda).$$
Put
$$B(\lambda)\buildrel \rm def \over =B_{1,1}(\lambda)=\ldots =
B_{1,d}(\lambda).$$
Estimates (2.12) follow from equalities (2.15), (2.16) and inequalities
(2.9), (2.10).

\vskip 2 mm
{\bf Proposition 2.3}

Let, under conditions of Proposition 2.2,  conductivity function
$\sigma\in C^{(3)}(V)$. Then\ $\forall\ z\in V$ function
$\lambda\mapsto\mu(z,\lambda)$, $\lambda\in\C$,is a unique solution of the
following integral equation
$$\mu(z,\lambda)=1-{1\over 2\pi i}\int\limits_{\xi\in\C}b(\xi)
e^{\bar\xi\bar z_1-\xi z_1}\overline{\mu(z,\xi)}{d\xi\wedge d\bar\xi\over
{\xi-\lambda}},\eqno(2.17)$$
where function $\lambda\mapsto b(\lambda)$ from (2.11) belongs to
$L^q(\C)$\ $\forall\ q\in (4/3,4)$.

\vskip 2 mm
{\it Proof}

Estimates (2.12) imply that\ $\forall\ q\in (4/3,4)$ function $b(\lambda)$
from Proposition 2.2 belongs to $L^q(\C)$. By the classical Vekua result
[Ve] with such functions $b(\lambda)$ the equation (2.17) is the uniquely
solvable Fredholm integral equation in the space $C(\bar{\C})$.

\vskip 2 mm
{\bf $\S 3$. Reconstruction of Faddeev type function $\psi$
and of conductivity $\sigma$ on $X$ through Dirichlet to
Neumann data on $bX$}

Let $X$ be domain with smooth (de class $C^{(2)}$) boundary on $V$ such that
$X\supseteq\bar V_0$.
Let $\sigma\in C^{(3)}(V)$, $\sigma>0$ on $V$ and $\sigma\equiv const$ on
$V\b X$. Let $q={d d^c\sqrt{\sigma}\over \sqrt{\sigma}}$, then
$q\in C_{1,1}^{(1)}(X)$ and $supp\,q\subseteq X$.

\vskip 2 mm
{\bf Definition} (Dirichlet-Neumann operator)

Let $u\in C^{(1)}(bX)$ and $U\in W^{1,\tilde p}(X)$,
$\tilde p>2$  be solution of the Dirichlet problem
$d\sigma d^cU\big|_X=0$, $U\big|_{bX}=u$,
where $d^c=i(\bar\pa-\pa)$. Operator
$u\big|_{bX}\to\sigma d^cU\big|_{bX}$ is called usually by
Dirichlet to Neumann operator.
Put $\tilde\psi=\sqrt{\sigma}U$ and $\psi=\sqrt{\sigma}u$. Then
$$d d^c\tilde\psi={d d^c\sqrt{\sigma}\over \sqrt{\sigma}}\tilde\psi=
q\tilde\psi\ \ {\rm on}\ \ X.\eqno(3.1)$$
Operator $\psi\big|_{bX}\mapsto\bar\pa\tilde\psi\big|_{bX}$ will be called
here by Dirichlet-Neumann operator.
The data contained in operator $u\big|_{bX}\to\sigma d^cU\big|_{bX}$
and in $\psi\big|_{bX}\to\bar\pa\tilde\psi\big|_{bX}$ are equivalent, but
for our further statements operator
$\psi\big|_{bX}\to\bar\pa\tilde\psi\big|_{bX}$ is more convenient.

Let $\psi_0$ be solution of Dirichlet problem
$$d d^c\psi_0=0,\ \ \psi_0\big|_{bX}=\psi\big|_{bX}.$$
Put
$$\hat\Phi\psi=\bar\pa\tilde\psi\big|_{bX}\ \ {\rm and}\ \
\hat\Phi_0\psi=\bar\pa\tilde\psi_0\big|_{bX},\eqno(3.2)$$
where operator $\psi\big|_{bX}\to\hat\Phi_0\psi$ is  Dirichlet-Neumann
operator for equation (3.1) with potential $q\equiv 0$.

\vskip 2 mm
{\bf Proposition 3.1} (Reconstruction of $\psi\big|_{bX}$ through
Dirichlet-Neumann data)

Let $\psi=\mu(z,\lambda)\cdot e^{\lambda z_1}$ be the Faddeev function
associated with potential $q={d d^c\sqrt{\sigma}\over \sqrt{\sigma}}$. Then
\ $\forall\lambda\in\C\b\{0\}$ the restriction
$\psi\big|_{bX}$ of $\psi$ on $bX$ is a unique solution in $C(bX)$
of the Fredholm integral equation:
$$\psi(z,\lambda)\big|_{bX}=e^{\lambda z_1}-
\int\limits_{\xi\in bX}e^{\lambda(z_1-\xi_1)}g_{\lambda}(z,\xi)\cdot
(\hat\Phi\psi(\xi)-\hat\Phi_0\psi(\xi)),\eqno(3.3)$$
where $g_{\lambda}(z,\xi)$ - kernel of operator $R_{\lambda}\circ\hat R$,
$$\hat\Phi\psi(\xi)-\hat\Phi_0\psi(\xi)=
\int\limits_{w\in bX}(\Phi(\xi,w)-\Phi_0(\xi,w))\psi(w),$$
$\Phi(\xi,w)$, $\Phi_0(\xi,w)$ are kernels of operators $\hat\Phi$ and
$\hat\Phi_0$.

\vskip 2 mm
{\it Remark}

This proposition for the case $V=\C$ coincide with the second part of
Theorem 1 from [N1].

\vskip 2 mm
{\bf Lemma 3.1} (Green-Riemann formula)

$\forall\ f,g\in C^{(1)}(X)$ we have equality
$$\int\limits_Xg\wedge\pa\bar\pa f-\int\limits_Xf\wedge\pa\bar\pa g=
\int\limits_{bX}g\wedge\bar\pa f+\int\limits_{bX}f\wedge\pa g.$$

\vskip 2 mm
{\it Proof}
$$\eqalign{
&\int\limits_Xg\wedge\pa\bar\pa f-\int\limits_Xf\wedge\pa\bar\pa g=
\int\limits_Xg\wedge\pa\bar\pa f+\int\limits_Xf\wedge\bar\pa\pa g=\cr
&\int\limits_{bX}g\wedge\bar\pa f-\int\limits_X\pa g\wedge\bar\pa f
+\int\limits_{bX}f\wedge\pa g-\int\limits_X\bar\pa f\wedge\pa g=\cr
&\int\limits_{bX}g\wedge\bar\pa f+\int\limits_{bX}f\wedge\pa g.\cr}$$

\vskip 2 mm
{\bf Lemma 3.2}

Let $\psi=e^{\lambda z_1}\mu$ be the Faddeev function associated
with $q\in C_{1,1}(V)$, $supp\,q\subseteq X$.
Let $G_{\lambda}(z,\xi)=e^{\lambda(z_1-\xi_1)}g_{\lambda}(z,\lambda)$, where
$g_{\lambda}(z,\xi)$ - kernel of operator $R_{\lambda}\circ\hat R$. Then

\noindent
$\forall\ z\in V\b X$ we have equality
$$\psi(z,\lambda)=e^{\lambda z_1}-\int\limits_{\xi\in bX}
G_{\lambda}(z,\xi)\bar\pa\psi(\xi)-
\int\limits_{\xi\in bX}\psi(\xi)\pa G_{\lambda}(z,\xi).\eqno(3.4)$$

\vskip 2 mm
{\it Proof}

From definition of the Faddeev function $\psi=e^{\lambda z_1}\mu$ we have
$$\psi(z,\lambda)=e^{\lambda z_1}+\hat G_{\lambda}\bigl({i\over 2}q\psi\bigr)
\ \ {\rm on}\ \ V,\eqno(3.5)$$
where $\hat G_{\lambda}$ - operator with kernel $G_{\lambda}(z,\xi)$
and $d d^c\psi=q\psi$. Besides,
$$\int\limits_XG_{\lambda}\bigl({i\over 2}q\psi\bigr)=
-\int\limits_XG_{\lambda}\pa\bar\pa\psi.$$
Using the Green-Riemann formula from Lemma 3.1 we have
$$\int\limits_XG_{\lambda}\pa\bar\pa\psi=
\int\limits_X\psi\pa\bar\pa G_{\lambda}+
\int\limits_{bX}G_{\lambda}\bar\pa\psi+
\int\limits_{bX}\psi\pa G_{\lambda}.$$
For $z\in V\b X$ we have $\pa\bar\pa G_{\lambda}=0$. Hence,
$$-\int\limits_XG_{\lambda}\bigl({i\over 2}q\psi\bigr)=
\int\limits_{bX}G_{\lambda}\bar\pa\psi+
\int\limits_{bX}\psi\pa G_{\lambda}.\eqno(3.6)$$
From (3.5) and (3.6) we obtain (3.4).

\vskip 2 mm
{\it Proof of Proposition 3.1}

Let $\psi_0:\ \bar\pa\pa\psi=0$  and $\psi_0\big|_{bX}=\psi$. Then by
Lemma 3.1\ $\forall\ z\in V\b X$ we have
$$\int\limits_{bX}\psi_0\pa G_{\lambda}+\int\limits_{bX}
G_{\lambda}\bar\pa\psi_0=0.\eqno(3.7)$$
Combining this relation with (3.4)   we obtain
$$\psi(z,\lambda)=e^{\lambda z_1}-\int\limits_{bX}G_{\lambda}(z,\xi)
(\bar\pa\psi(\xi)-\bar\pa\psi_0(\xi)).\eqno(3.8)$$
Equality (3.8) implies (3.3), where $\hat\Phi\psi$, $\hat\Phi_0\psi$ are
defined by (3.2). Integral equation (3.3) is the Fredholm equation in
$L^{\infty}(bX)$, because operator $(\hat\Phi-\hat\Phi_0)$ is a compact
operator in
$L^{\infty}(bX)$. For the case $V=\C$ Proposition 1 of [N1], contains
explicit inequality:
$$\lim\limits_{\scriptstyle w\to\xi \atop \scriptstyle w,\xi\in bX}
{|\Phi(\xi,w)-\Phi_0(\xi,w)|\over |\ln|\xi-w||}\le
(const)\big|{d d^c\sqrt{\sigma}\over \sqrt{\sigma}}\big|(w).\eqno(3.9)$$
Because
of its local nature, equality (3.9) is valid for our more general case.

Existence of the Faddeev function\ $\forall\lambda\ne 0$,
proved in Proposition 1.2, implies existence of solution of
equation (3.3). The uniqueness of solution of equation (3.3)
we deduce (as in  proof of Proposition 2 in [N1]) from the  following
statement.

\vskip 2 mm
{\bf Lemma 3.3}

Let $q\in C_{1,1}(V)$, $supp\,q\subseteq X$. Then each solution
$\psi=\psi(z,\lambda)$ of equation $d d^c\psi=q\psi$ on $V$ which coincides
on $bX$ with solution of integral equation (3.3) in space $C(bX)$ is the
Faddeev type function associated with $q$.

\vskip 2 mm
{\it Proof}

Let $\psi$ satisfy $d d^c\psi=q\psi$ on $V$, $\psi\big|_{bX}\in C(bX)$ and
$\psi\big|_{bX}$ satisfy integral equation (3.3)=(3.8). From (3.8) with
help of (3.7) we obtain formula (3.4) for $\psi(z,\lambda)$,
$z\in V\b X$, $\lambda\in\C$. From the Sohotsky-Plemelj jump formula
we deduce (see [HM], Lemma 15) that\ $\forall\ z^*\in bX$
$$\psi(z^*)=\lim\limits_{\scriptstyle z\to z^* \atop \scriptstyle z\in X}
\bigl(\int\limits_{bX}G_{\lambda}\bar\pa\psi+\psi\pa G_{\lambda}\bigr)-
\lim\limits_{\scriptstyle z\to z^* \atop \scriptstyle z\in V\b X}
\bigl(\int\limits_{bX}G_{\lambda}\bar\pa\psi+\psi\pa G_{\lambda}\bigr).
\eqno(3.10)$$
From (3.4) and (3.10) we obtain equality
$$e^{\lambda z_1}=\int\limits_{bX}G_{\lambda}\bar\pa\psi+\psi\pa G_{\lambda}\
\ {\rm for}\ \ z\in X,\ \ \lambda\in\C.\eqno(3.11)$$
Using the Green-Riemann formula (Lemma 3.1) we obtain further
$$\eqalign{
&\int\limits_{bX}G_{\lambda}\bar\pa\psi+\psi\pa G_{\lambda}=
\int\limits_X\psi\bar\pa\pa G_{\lambda}-
\int\limits_XG_{\lambda}\bar\pa\pa\psi=\cr
&\left\{\matrix{
\psi(z)-\int\limits_XG_{\lambda}\bar\pa\pa\psi,&{\rm if}\ z\in X,\
\lambda\in\C \hfill\cr
-\int\limits_XG_{\lambda}\bar\pa\pa\psi,&{\rm if}\ z\in V\b X,\
\lambda\in\C \hfill\cr}\right.\cr}\eqno(3.12)$$
Equalities (3.4), (3.11) and (3.12) imply
$$\psi(z)=e^{\lambda z_1}+\int\limits_VG_{\lambda}\bar\pa\pa\psi=
e^{\lambda z_1}+\hat G_{\lambda}({i\over 2}q\psi),$$
i.e. $\psi$ satisfies (3.5).
From equality (3.5) and Proposition 3 ii) from [H2] we obtain that
$\mu-1=\psi e^{-\lambda z_1}-1\in W^{1,\tilde p}(V)$ and $\psi$
is the Faddeev type function associated with $q$.

Lemma 3.3 is proved.

The uniqueness of solution of equation (3.3) in $C(bX)$ follows now from
Lemma 3.3 and the uniqueness of the Faddeev type function associated with
$q$ proved in

\noindent
Proposition 1.1.

Proposition 3.1 is proved.

\vskip 2 mm
{\bf Proposition 3.2} (Reconstruction of
$\psi\big|_X$ through $\psi\big|_{bX}$)

In conditions of Propositions 1.1, 1.2 the following properties for
"$\bar\pa$-scattering data $b(\lambda)$" permit to reconstruct
$\psi\big|_X$ through $\psi\big|_{bX}$
$$
\forall\ z^*\in bX\ \exists\lim\limits_{\scriptstyle z\to\infty \atop
\scriptstyle z\in V} {\bar z_1\over \bar\lambda}e^{-\bar\lambda\bar z_1}
{\pa\psi\over \pa\bar z_1}(z,\lambda)\buildrel \rm def \over =b(\lambda)=
\bar\psi(z^*,\lambda)^{-1}{\pa\psi\over \pa\bar\lambda}(z^*,\lambda),$$
$\forall\ z\in X$ function $\lambda\mapsto\psi(z,\lambda)$, $\lambda\in\C$,
is a unique solution of the integral equation
$$\psi(z,\lambda)=e^{\lambda z_1}-{1\over 2\pi i}\int\limits_{\xi\in\C}
b(\xi)e^{(\lambda-\xi)z_1}\bar\psi(\xi,\lambda){d\xi\wedge d\bar\xi\over
{\xi-\lambda}}.$$
This proposition is a consequence of Propositions 2.1-2.3.

\vskip 2 mm
{\bf Proposition 3.3} (Reconstruction formulas for $\sigma$)

Conductivity function $\sigma\big|_X$ can be reconstructed through
Dirichlet-Neumann data by the R.Novikov's scheme:
$${\rm DN\ data}\ \to\psi\big|_{bX}
\to\bar\pa-{\rm scattering\ data}
\to\psi\big|_X\to
{d d^c\sqrt{\sigma}\over \sqrt{\sigma}}\big|_X.$$
The last step of this scheme one can realize by any of the following
formulas:
$$\eqalign{
&{\rm A})\ {d d^c\sqrt{\sigma}\over \sqrt{\sigma}}(z)=
(d d^c\psi(z,\lambda))\psi^{-1}(z,\lambda),\ \ z\in X,\ \lambda\in\C,\cr
&{\rm B})\ {d d^c\sqrt{\sigma}\over \sqrt{\sigma}}(z)=
2i\lim\limits_{\lambda\to\infty}\lambda e^{-\lambda z_1}dz_1\wedge
\bar\pa\psi(z,\lambda),\ z\in X,\cr
&{\rm C})\ {d d^c\sqrt{\sigma}\over \sqrt{\sigma}}(z)=
-2i\lim\limits_{\lambda\to 0}{\bar\pa\pa\mu(z,\lambda)\over \mu(z,\lambda)},\
 z\in X.\cr}$$

\vskip 2 mm
{\it Proof}

The first three steps of the scheme above were done in Propositions 3.1, 3.2.
The formula A) is an immediate consequense of equation (3.1).
The formula B) follows from equation
$\bar\pa(\pa+\lambda dz_1)\mu={i\over 2}
{d d^c\sqrt{\sigma}\over \sqrt{\sigma}}\mu$ (Proposition 1.2a) and
estimates  $\mu\to 1$, $\lambda\to\infty$

\noindent
(Proposition 1.2b),
$\pa\bar\pa\mu(z,\lambda)\to 0$, $\lambda\to\infty$ (Proposition 1.2e).
The formula C)  follows from the same equation and estimates
$\|{\pa\mu\over \pa z_1}\|=O(\sqrt{|\lambda|}),\  \lambda\to 0$,\
(Proposition 1.2c).

\vskip 4 mm
{\bf References}

\item{[ BC1]} Beals R., Coifman R., Multidimensional inverse scattering and
nonlinear partial differential equations, Proc.Symp. Pure Math. {\bf 43}
(1985), A.M.S. Providence, Rhode Island, 45-70
\item{[ BC2]} Beals R., Coifman R., The spectral problem for the
Davey-Stewartson  and Ishimori hierarchies, In: "Nonlinear Evolution Equations
: Integrability and Spectral Methodes", Proc. Workshop, Como, Italy 1988,
Proc. Nonlinear Sci., 15-23, 1990
\item{[  Be]} Belishev M.I., The Calderon problem for two dimensional
manifolds by the BC-method, SIAM J.Math.Anal. {\bf 35}:1, (2003), 172-182
\item{[  BU]} Brown R., Uhlmann G., Uniqueness in the inverse conductivity
problem for non-smooth conductivities in two dimensions, Comn.Part.Dif.
Equations {\bf 22} (1997), 1009-1027
\item{[   C]} Calderon A.P., On an inverse boundary problem, Seminar on
Numerical Analysis and its Applications to Continuum Physics, Soc.
Brasiliera de Matematica, Rio de Janeiro (1980), 61-73
\item{[  F1]} Faddeev L.D., Increasing solutions of the Schr\"odinger
equation, Dokl.Akad.Nauk SSSR {\bf 165} (1965), 514-517 (in Russian);
Sov.Phys.Dokl. {\bf 10} (1966), 1033-1035
\item{[  F2]} Faddeev L.D., The inverse problem in the quantum theory of
scattering, II, Current Problems in Math., {\bf 3}, 93-180, VINITI, Moscow,
1974 (in Russian); J.Soviet Math. {\bf 5} (1976), 334-396
\item{[   G]} Gel'fand I.M., Some problems of functional analysis and
algebra, Proc.Internat.
\item{      } Congr.Math., Amsterdam (1954), 253-276
\item{[  GN]} Grinevich P.G., Novikov S.P., Two-dimensional "inverse
scattering problem" for negative energies and generalized analytic functions,
Funct.Anal. and Appl., {\bf 22} (1988), 19-27
\item{[  H1]} Henkin G.M., Two-dimensional electrical tomography and
complex analysis, Lecture at the Mathematical weekend, European Math. Soc.,
Nantes, june 16-18, 2006,
http://univ-nantes.fr/WEM2006
\item{[  H2]} Henkin G.M., Cauchy-Pompeiu type formulas for $\bar\pa$
on affine algebraic Riemann surfaces and some applications, Preprint 2008,
arXiv:0804.3761

\item{[  HM]} Henkin G.M., Michel V., On the explicit reconstruction of a
Riemann surface from its Dirichlet-Neumann operator, {\bf 17} (2007), 116-155
\item{[  LU]} Lassas M., Uhlmann G., On determining a Riemann manifold from
the Dirichlet-to-Neumann map, Ann.Sci.Ecole Norm.Sup. {\bf 34} (2001),
771-787
\item{[  N1]} Novikov R., Multidimensional inverse spectral problem for the
equation

\item{      }
$-\Delta\psi+(v-Eu)\psi=0$, Funkt.Anal. i Pril. {\bf 22} (1988),
11-22 (in Russian); Funct.Anal and Appl. {\bf 22} (1988), 263-278
\item{[  N2]} Novikov R., The inverse scattering problem on a fixed energy
level for the two-dimensional Schr\"odinger operator, J.Funct.Anal. {\bf 103}
(1992), 409-463
\item{[  N3]} Novikov R., Reconstruction of a two-dimensional Schr\"odinger
operator from the scattering amplitude at fixed energy,
Funkt.Anal. i Pril. {\bf 20} (1986),
90-91 (in Russian); Funct.Anal and Appl. {\bf 20} (1986), 246-248
\item{[  Na]} Nachman A., Global uniqueness for a two-dimensional inverse
boundary problem, Ann. of Math. {\bf 143} (1996), 71-96
\item{[   R]} Rodin Yu., Generalized analytic functions on Riemann surfaces,
Lecture Notes Math., 1288, Springer 1987
\item{[   S]} Sylvester J., An anisotropic inverse boundary value problem,
Comm.Pure Appl.Math. {\bf 43} (1990), 201-232
\item{[  SU]} Sylvester J., Uhlmann G., A uniqueness theorem for an inverse
boundary value problem in electrical prospection, Comm.Pure Appl.Math.{\bf 39}
(1986), 91-112
\item{[   T]} Tsai T.Y., The Shr\"odinger equation in the plane, Inverse
problems, {\bf 9} (1993), 763-787
\item{[   V]} Vekua I.N., Generalized analytic functions, Pergamon Press, 1962

\vskip 4 mm
\line{\hfill \box3}

\end